\documentclass[10pt,a4paper,reqno]{amsart} 
\usepackage{pstricks,amsmath,here,latexsym,amssymb}
\usepackage{epsf,psfrag,psfig,epsfig}
\usepackage[latin1]{inputenc}
\usepackage[english]{babel} 
\usepackage{t1enc}

\newtheorem{thm}{Theorem}[section]
\newtheorem{prop}[thm]{Proposition}
\newtheorem{cor}[thm]{Corollary}
\newtheorem{lemma}[thm]{Lemma}
\newcommand{\dem}{\noindent \textbf{Proof: }}

\theoremstyle{definition} 
\newtheorem{Def}[thm]{Definition}

\newcommand{\Ref}[1]{(\ref{#1})}
\newcommand{\findem}{\vspace{-.4cm} \begin{flushright} $\square~$ \end{flushright} \vspace{.4cm} }

\newcommand{\titre}[1]{\noindent \textbf{#1}}

\newcommand{\ite}{\noindent $\bullet~$}
\newcommand{\iten}{\noindent -~}

\author[O. Bernardi]{Olivier Bernardi}
\address{LaBRI, Université Bordeaux 1, 351 cours de la Libération, 33405 Talence Cedex, France} 

\email{bernardi@labri.fr}

\title[Bijective counting of Kreweras walks and loopless triangulations]{Bijective counting of Kreweras walks and loopless triangulations}

\subjclass[2000]{Primary 05A15}
\keywords{planar walk, Kreweras walk, planar map, triangulation, cubic map, bijection, counting}

\begin{document}

\begin{abstract}
We consider lattice walks in the plane starting at the origin, remaining in the first quadrant $i,j\geq 0$ and made of West, South and North-East steps. In 1965, Germain Kreweras discovered a remarkably simple formula giving the number of these walks (with prescribed length and endpoint). Kreweras' proof was very involved and several alternative derivations have been proposed since then. But the elegant simplicity of the counting formula remained unexplained. We give the first purely combinatorial explanation of this formula. Our approach is based on a bijection between Kreweras walks and triangulations with a distinguished spanning tree. We obtain simultaneously a bijective way of counting loopless triangulations.
\end{abstract}

\maketitle

\section{\textbf{Introduction}}
\begin{figure}[htb!]
\begin{center}
\input{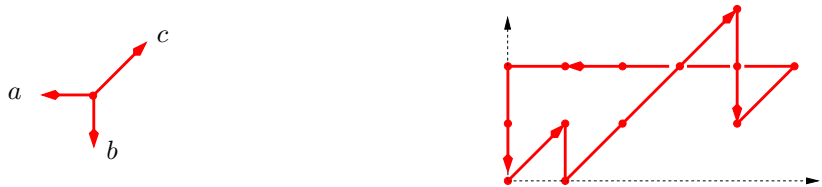}
\caption{The Kreweras walk $cbcccbbcaaaaabb$.}\label{fig:walk}
\end{center}
\end{figure}

We consider lattice walks in the plane starting from the origin (0,0), remaining in the first quadrant $i,j\geq 0$ and made of three kind of steps: West, South and North-East. These walks were first studied by Kreweras \cite{Kreweras:walks} and inherited his name. A Kreweras walk ending at the origin is represented in Figure~\ref{fig:walk}. 

These walks have remarkable enumerative properties. Kreweras proved in 1965 that the number of walks of length $3n$ ending at the origin is:
\begin{eqnarray}
k_n &=& \frac{4^n}{(n+1)(2n+1)} {3n \choose n}. \label{eq:Kreweras}
\end{eqnarray}
The original proof of this result is complicated and somewhat unsatisfactory. It was performed by \emph{guessing} the number of walks of size $n$ ending at a generic point $(i,j)$. The conjectured formulas were then checked using the recurrence relations between these numbers. The \emph{checking part} involved several hypergeometric identities which were later simplified by Niederhausen \cite{Niederhausen:Kreweras}.  In 1986, Gessel gave a different proof in which the \emph{guessing part} was reduced \cite{Gessel:Kreweras}. 
More recently,  Bousquet-Mélou proposed a constructive proof (that is, without guessing) of these results and some extensions \cite{MBM:Kreweras}. Still, the simple looking formula \Ref{eq:Kreweras} remained without a direct combinatorial explanation. The problem of finding a combinatorial explanation was mentioned by Stanley in \cite{Stanley-Clay-problems}. Our main goal in this paper is to provide such an explanation. \\

Formula \Ref{eq:Kreweras} for the number of Kreweras walks is to be compared to another formula proved the same year. In 1965, Mullin, following the seminal steps of Tutte, proved  via a generating function approach \cite{Mullin:triangulation-nonsep} that the number of loopless triangulations of size $n$ (see below for precise definitions) is 
\begin{eqnarray}
t_n &=& \frac{2^n}{(n+1)(2n+1)} {3n \choose n}. \label{eq:triangulations}
\end{eqnarray}
A bijective proof of \Ref{eq:triangulations} was outlined by Schaeffer in his Ph.D thesis \cite{Schaeffer:these}. See also \cite{Schaeffer:triangulation} for a more general construction concerning  loopless triangulations of a $k$-gon. We will give an alternative bijective proof for the number of loopless triangulations. Technically speaking, we will work instead on \emph{bridgeless cubic maps} which are the dual of loopless  triangulations.\\

It is interesting to observe that both \Ref{eq:Kreweras} and \Ref{eq:triangulations} admit a nice generalization. Indeed, the number $k_{n,i}$ of Kreweras walks of size $n$ ending at point $(i,0)$ and  the number $c_{n,i}$ of loopless triangulations of size $n$ of an $(i+2)$-gon both admit a closed  formula (see \Ref{eq:kni} and \Ref{eq:tni}). Moreover, the numbers $k_{n,i}$ and $c_{n,i}$ are related by the equation $k_{n,i}=2^n c_{n,i}$. This relation is explained in Section \ref{section:perspectives}. Alas, we have found no way of proving these formulas by our approach. \\

\section{\textbf{How the proofs work}}
\begin{figure}[h!]
\begin{center}
\input{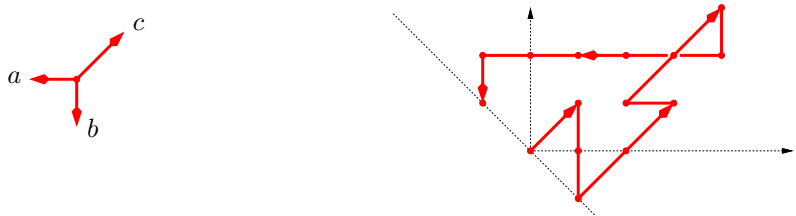}
\caption{An excursion.}\label{fig:extended-walk}
\end{center}
\end{figure}

We begin with an account of this paper's content in order to underline the (slightly unusual) logical structure of our proofs. \vspace{.2cm}\\
\ite In Section \ref{section:preliminaries}, we first recall some definitions on planar maps. We also define a special class of spanning trees called \emph{depth trees}. Depth trees are closely related to the trees that can be obtained by a \emph{depth first search} algorithm. \vspace{.2cm}\\
Then, we consider a larger family of walks containing the Kreweras walks. These walks are made of West, South and North-East steps, start from the origin and remain in the half-plane $i+j\geq 0$. We borrow a terminology from probability theory and call these walks \emph{meanders}. We call \emph{excursion} a meander ending on the second diagonal (i.e. the line $i+j=0$). An excursion is represented in Figure~\ref{fig:extended-walk}.\\

\noindent Unlike Kreweras walks, excursions are easy to count. By applying the cycle lemma (see \cite[Section 5.3]{Stanley:volume2}), we prove that the number of excursions of size $n$ (length $3n$) is 
$$e_n=\frac{4^n}{2n+1}{3n \choose n}~.$$
\ite In Section \ref{section:extension-bijection}, we define a mapping   $\Phi$ between excursions and cubic maps with a distinguished depth tree. In Section \ref{section:preuve-bijection} we prove that the mapping $\Phi$ is a $(n+1)$-to-1 correspondence $\Phi$ between excursions (of size $n$) and bridgeless cubic maps (of size $n$) with a distinguished depth tree. As a consequence, the number of bridgeless cubic maps of size $n$ with a distinguished depth tree is found to be:
$$d_n~=~\frac{e_n}{n+1}~=~\frac{4^n}{(n+1)(2n+1)}{3n \choose n}~.$$
\ite In Section \ref{section:bijection}, we prove that the correspondence $\Phi$, restricted to Kreweras walks, induces a \emph{bijection} between Kreweras walks (of size $n$)  ending at the origin and bridgeless cubic maps (of size $n$) with a distinguished depth tree. As a consequence, we obtain:
$$k_n~=~d_n~=~\frac{4^n}{(n+1)(2n+1)}{3n \choose n},$$
where $k_n$ is the number of Kreweras walks of size $n$ ending at the origin. This gives a combinatorial proof of \Ref{eq:Kreweras}.\vspace{.2cm}\\
\ite In Section \ref{section:counting-cubic}, we enumerate depth trees on cubic maps. We prove that the number of such trees for a cubic map of size $n$ is $2^n$. As a consequence, the number of cubic maps of size $n$ is 
$$c_n~=~\frac{d_n}{2^n}~=~\frac{2^n}{(n+1)(2n+1)}{3n \choose n}.$$
This gives a combinatorial proof of \Ref{eq:triangulations}. \vspace{.2cm}\\
\ite In Section \ref{section:perspectives}, we extend the mapping $\Phi$ to Kreweras walks ending at $(i,0)$ and discuss some open problems.\\

\section{\textbf{Preliminaries}}  \label{section:preliminaries}
\subsection{Planar maps and depth trees} \label{section:maps}
~\\
\titre{Planar maps.}
A \emph{planar map}, or \emph{map} for short, is an embedding of a connected planar graph in the sphere without intersecting edges,  defined up to orientation preserving homeomorphisms of the sphere. Loops and multiple edges are allowed. The \emph{faces} are the connected components of the complement of the graph. By removing the midpoint of an edge we obtain two \emph{half-edges}, that is, one-dimensional cells incident to one vertex. We say that each edge has two \emph{half-edges}, each of them incident to one of the endpoints. \\

\begin{figure}[h!]
\begin{center}
\input{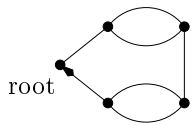}
\caption{A rooted map.}\label{fig:example-map}
\end{center}
\end{figure}

A map is \emph{rooted} if one of its half-edges is distinguished as the \emph{root}. The edge containing the root is the \emph{root-edge} and its endpoint is the \emph{root-vertex}. Graphically, the root is indicated by an arrow pointing on the root-vertex (see Figure \ref{fig:example-map}). All the maps considered in this paper are rooted and we shall not further precise it. \\

\begin{figure}[h!]
\begin{center}
\input{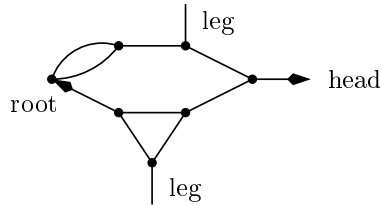}
\caption{A growing map.}\label{fig:example-growing-map}
\end{center}
\end{figure}

\titre{Growing maps.}
Our constructions lead us to consider maps with some \emph{legs}, that is, half-edges that are not part of a complete edge. A \emph{growing map} is a (rooted) map together with some legs, one of them being distinguished as  the \emph{head}. We require the legs to be (all) in the same face called \emph{head-face}. The endpoint of the head is the \emph{head-vertex}. Graphically, the head is indicated by an arrow pointing away from the head-vertex. The root of a growing map can be a leg or a regular half-edge. For instance, the growing map in Figure~\ref{fig:example-growing-map} has 2 legs beside the head, and its root is not a leg.\\

\begin{figure}[h!]
\begin{center}
\input{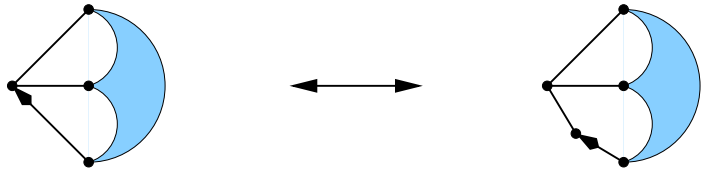}
\caption{Bijection between cubic maps and 2-near-cubic maps.}\label{fig:2-near-cubic}
\end{center}
\end{figure}

\titre{Cubic maps.}
A  map (or growing map) is \emph{cubic} if every vertex has degree 3. It is $k$-near-cubic if the root-vertex has degree $k$ and any other vertex has degree 3. For instance, the  map in Figure~\ref{fig:example-map} is 2-near-cubic and the  growing map in  Figure~\ref{fig:example-growing-map} is cubic. Observe that cubic maps are in bijection with 2-near-cubic maps not reduced to a loop by the mapping illustrated in Figure \ref{fig:2-near-cubic}.\\

The incidence relation between vertices and edges in cubic maps shows that the number of edges is always a multiple of $3$. More generally, if $M$ is a $k$-near-cubic map with $e$ edges and $v$ vertices, the incidence relation reads: $3(v-1)+k=2e$. Equivalently, $3(v-k+1)=2(e-2k+3)$. The number $v-k+1$ is non-negative for \emph{non-separable} $k$-near-cubic maps (see definition below). (This property can be shown by induction on the number of edges by contracting the root-edge.) Hence, the number of edges has the form $e=3n+2k-3$, where $n$ is a non-negative integer. We say that  a $k$-near-cubic map has \emph{size} $n$ if it has  $e=3n+2k-3$ edges (and $v=2n+k-1$ vertices). In particular, the mapping of Figure \ref{fig:2-near-cubic} is a bijection between cubic maps of size $n$ ($3n+3$ edges) and 2-near-cubic maps of size $n+1$ ($3n+4$ edges).\\

\titre{Non-separable maps.}
A map is \emph{non-separable} if its edge set cannot be partitioned into two non-empty parts such that only one vertex is incident to some edges in both parts. In particular, a non-separable map not reduced to an edge has no loop nor bridge (a \emph{bridge} or \emph{isthmus} is an edge whose deletion  disconnects the map). For cubic maps and 2-near-cubic maps it is equivalent to be non-separable or bridgeless.  The mapping illustrated in Figure \ref{fig:2-near-cubic} establishes a bijection between bridgeless cubic maps and bridgeless 2-near-cubic maps not reduced to a loop.\\

\begin{figure}[h!]
\begin{center}
\input{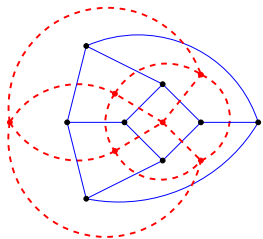}
\caption{A cubic map and the dual triangulation (dashed lines).}\label{fig:cubic-triangulation}
\end{center}
\end{figure}

Bridgeless cubic maps are interesting because their \emph{dual} are the loopless triangulations. Recall that the dual $M^*$ of a map $M$ is the map obtained by putting a vertex of $M^*$ in each face of $M$ and an edge of $M^*$ across each edge of $M$. See Figure~\ref{fig:cubic-triangulation} for an example.\\

\begin{figure}[h!]
\begin{center}
\input{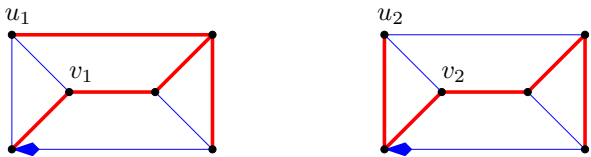}
\caption{A depth tree (left) and a non-depth tree (right).}\label{fig:depth-three}
\end{center}
\end{figure}

\titre{Depth trees.}
A \emph{tree} is a connected graph without cycle. A subgraph $T$ of a graph $G$ is a \emph{spanning tree} if it is a tree containing every vertex of $G$. An edge of the graph $G$ is said to be \emph{internal} if it is in the spanning tree $T$ and \emph{external} otherwise. For any pair of vertices $u,~v$ of the graph $G$, there is a unique path between $u$ and $v$ in the spanning tree $T$. We call it the \emph{$T$-path} between $u$ and $v$. A map (or growing map) $M$ with a distinguished spanning tree $T$ will be denoted by $M_T$. Graphically, we shall indicate the spanning tree  by thick lines as in Figure \ref{fig:depth-three}. A vertex $u$ of $M_T$ is an \emph{ancestor} of another vertex $v$ if it is on the $T$-path  between the root-vertex and $v$. In this case, $v$ is a \emph{descendant} of $u$. Two vertices are \emph{comparable} if one is the ancestor of the other.  For instance, in Figure \ref{fig:depth-three}, the vertices $u_1$ and  $v_1$ are comparable whereas $u_2$ and $v_2$ are not. \\

A \emph{depth tree} is a spanning tree such that any external edge joins comparable vertices. Moreover, we require the edge containing the root to be external. In Figure \ref{fig:depth-three}, the tree on the left side is a depth tree but the tree on the right side is not a depth tree since the edge $(u_2,v_2)$ breaks the rule. A \emph{depth-map} is a map with a distinguished depth tree. A \emph{marked-depth-map} is a depth-map with a marked external edge.\\

\subsection{Kreweras walks and meanders} \label{section:Kreweras-extended}~\\
In what follows, Kreweras walks are considered as words on the alphabet $\{a,b,c\}$. The letter $a$ (resp. $b$, $c$) corresponds to a West (resp. South, North-East) step. For instance, the walk in Figure \ref{fig:walk} is $cbcccbbcaaaaabb$. The length of a word $w$ is denoted by $|w|$ and the number of occurrences of a given letter $\alpha$ is denoted by $|w|_{\alpha}$. Kreweras walks are the words $w$ on the alphabet $\{a,b,c\}$ such that any prefix $w'$ of $w$ satisfies 
\begin{eqnarray} \label{eq:prefix}
|w'|_a \leq |w'|_c \hspace{.5cm} \textrm{and} \hspace{.5cm} |w'|_b \leq |w'|_c~. 
\end{eqnarray}
Kreweras walks ending at the origin satisfy the additional constraint 
\begin{eqnarray} \label{eq:word}
|w|_a~=~|w|_b~=~|w|_c. 
\end{eqnarray}
These conditions can be interpreted as a ballot problem with three candidates. This is why Kreweras walks sometimes appear under this formulation in the literature~\cite{Niederhausen:Kreweras}. \\

Similarly, the meanders, that is, the walks remaining in the half-plane $i+j\geq 0$, are the words $w$ on $\{a,b,c\}$ such that any prefix $w'$ of $w$ satisfies 
\begin{eqnarray} \label{eq:prefix-extended}
|w'|_a~ + ~|w'|_b ~\leq ~2 |w'|_c~.  
\end{eqnarray}
Excursions, that is, meanders ending on the second diagonal, satisfy the additional constraint 
\begin{eqnarray} \label{eq:word-extended}
|w|_a~+~|w|_b~=~2|w|_c~. 
\end{eqnarray}

Note that the length of any walk ending on the second diagonal is a multiple of~$3$. The \emph{size} of such a walk of length $3n$ is $n$. Note also that a walk ending at point $(i,0)$ has a length of the form $l=3n+2i$ where $n$ is a non-negative integer. A Kreweras walk of length  $l=3n+2i$ ending at $(i,0)$  has \emph{size} $n$.  \\

Unlike Kreweras walks, the  excursions are easy to count.

\begin{prop} \label{thm:count-extended}
There are 
\begin{eqnarray} e_n~=~\frac{4^n}{2n+1}{3n \choose n} \end{eqnarray}
excursions of size $n$.
\end{prop}

\begin{figure}[h!]
\begin{center}
\input{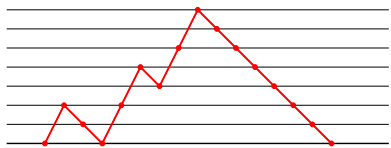}
\caption{The projected walk associated to the excursion of Figure~\ref{fig:extended-walk}.
}\label{fig:one-dimensionnal-walk}
\end{center}
\end{figure}

\dem
We consider \emph{projected walks}, that is, one-dimensional lattice walks starting and ending at $0$, remaining non-negative and made of steps $+2$ and $-1$.  (They correspond to projections of excursions on the first diagonal.) 
A projected walk is represented in Figure \ref{fig:one-dimensionnal-walk}. Projected walks can be seen as words $w$ on the alphabet $\{\alpha,c\}$ with  $|w|_\alpha = 2 |w|_c$ and such that any prefix $w'$ of $w$ satisfies $|w'|_\alpha \leq 2 |w'|_c$. 
The projected walks can be counted bijectively by applying the cycle lemma (see Section 5.3 of \cite{Stanley:volume2}): there are $$\displaystyle p_n=\frac{1}{3n+1}{3n+1 \choose 2n+1}=\frac{1}{2n+1}{3n \choose n}$$ 
projected walks of size $n$ (length $3n$).\\
Given an excursion, we obtain a projected walk by replacing the occurrences of $a$ and $b$ by $\alpha$. Conversely, taking a projected walk of length $3n$ and replacing the $2n$ letters $\alpha$ by a sequence of letters in $\{a,b\}$ one obtains an excursion. This establishes a $4^n$-to-1 correspondence between excursions (of size $n$)  and projected walks (of size $n$). Thus, there are $4^n p_n$ excursions of size $n$.
\findem

\section{\textbf{A bijection between excursions and cubic marked-depth-maps}}\label{section:extension-bijection}
\begin{figure}[h!]
\begin{center}
\input{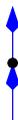}
\caption{The tree-growing map $M_\bullet^0$.}\label{fig:seed}
\end{center}
\end{figure}

\begin{figure}[h!]
\begin{center}
\input{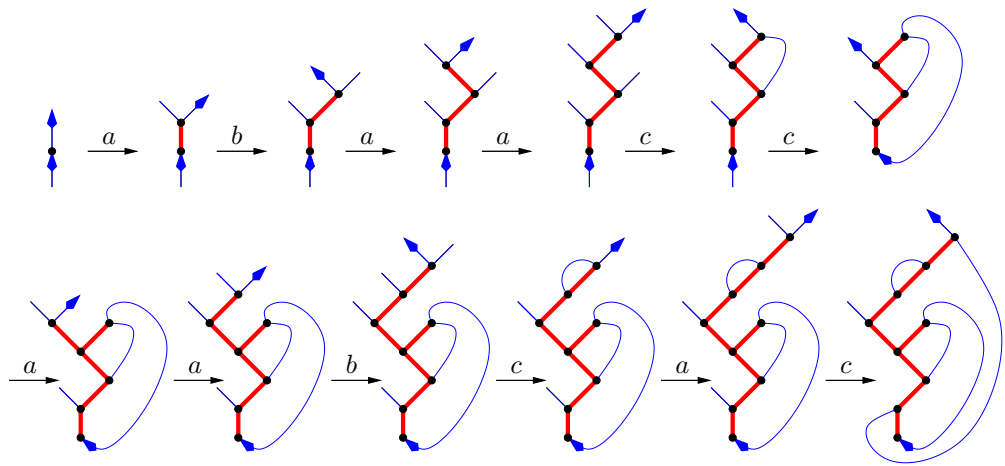}
\caption{Successive applications of the mappings $\varphi_a,~\varphi_b,~\varphi_c$ for the walk   $cacbaaccaaba$ (read from right to left).}\label{fig:example-extended-bijection}
\end{center}
\end{figure}

\begin{figure}[h!]
\begin{center}
\input{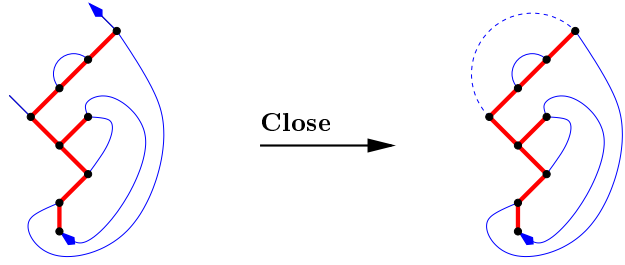}
\caption{Closing the map (the marked edge is dashed).}\label{fig:final-step}
\end{center}
\end{figure}

In this section we define a mapping $\Phi$ between excursions and bridgeless 2-near-cubic \emph{marked-depth-maps} (2-near-cubic maps with a distinguished depth tree and a marked external edge). We shall prove in Section \ref{section:preuve-bijection} that the mapping $\Phi$ is a bijection between  excursions and bridgeless 2-near-cubic marked-depth-maps. The general principle of the mapping $\Phi$ is to read the excursion from right to left and interpret each letter as an operation for constructing the map and the tree. This step-by-step construction is illustrated in Figure \ref{fig:example-extended-bijection}. The intermediate steps are \emph{tree-growing maps}, that is, growing maps together with a distinguished spanning tree (indicated by thick lines). \vspace{.2cm} \\
\ite We start with the tree-growing map $M_\bullet^0$ consisting of one vertex and two legs. One of the legs is the root, the other is the head (see Figure \ref{fig:seed}). The spanning tree is reduced to the unique vertex. \vspace{.2cm} \\
\ite We apply successively certain elementary mappings $\varphi_a,~\varphi_b,~\varphi_c$ (Definition \ref{def:varphi}) corresponding to the letters $a,b,c$ of the excursion read \emph{from right~to~left}.\vspace{.2cm} \\
\ite When the whole excursion  is read, there is only one leg remaining beside the head. At this stage, we \emph{close} the tree-growing map, that is, we glue the head and the remaining leg into a \emph{marked} external edge as shown in Figure \ref{fig:final-step}. \\  

\begin{figure}[h!]
\begin{center}
\input{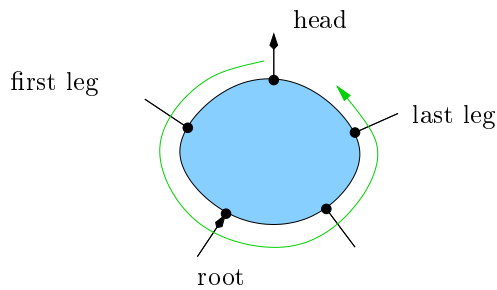}
\caption{Making the tour of the head-face.}\label{fig:making-the-tour}
\end{center}
\end{figure}

Let us enter in the details and define the mapping $\Phi$. Consider a growing map~$M$. We \emph{make a tour} of the head-face if we follow its border in counterclockwise direction (i.e. the border of the head-face stays on our left-hand side) starting from the head (see Figure \ref{fig:making-the-tour}). This journey induces a linear order on the legs of $M$.  We shall talk about the \emph{first} and \emph{last} legs of $M$. \\

\begin{figure}[h!]
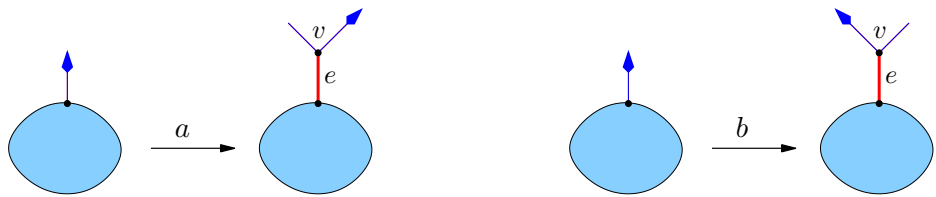

\begin{center}
\input{varphi-a.pstex_t}\hspace{2.5cm} \input{varphi-b.pstex_t}
\caption{The mappings $\varphi_a$ and $\varphi_b$.}\label{fig:varphi-a}
\end{center}
\end{figure}

\begin{figure}[h!]
\begin{center}
\input{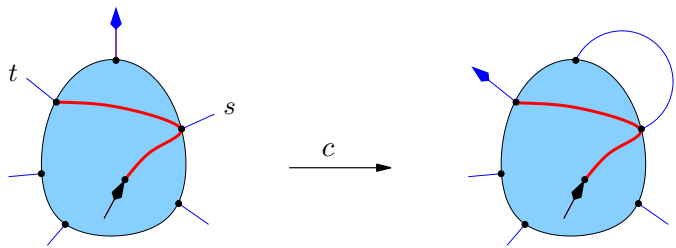}
\caption{The mapping $\varphi_c$.}\label{fig:varphi-c}
\end{center}
\end{figure}

We define three mappings $\varphi_a$,  $\varphi_b$, $\varphi_c$ on tree-growing maps.
\begin{Def} \label{def:varphi}
Let $M_T$ be a tree-growing map (the map is $M$ and the distinguished tree is $T$). \vspace{.2cm}\\
\ite The mappings $\varphi_a$ and $\varphi_b$ are represented in Figure \ref{fig:varphi-a}. The tree-growing map $M'_{T'}=\varphi_a(M_T)$ (resp. $\varphi_b(M_T)$) is obtained from $M_T$ by replacing the head by an edge $e$ together with a new vertex $v$ incident to the new head and another leg at its left (resp. right). The tree $T'$ is obtained from $T$ by adding the edge $e$ and the vertex $v$.\vspace{.2cm} \\
\ite The tree-growing map $\varphi_c(M_T)$ is only defined if the first and last legs exist (that is, if the head-face contains some legs beside the head) and have distinct and comparable endpoints. We call these legs $s$ and $t$ \emph{with the convention that the endpoint of $s$ is an ancestor of the endpoint of $t$}. \\
In this case, the tree-growing map $M'_{T}=\varphi_c(M_T)$ is obtained from $M_T$ by gluing together the head and the leg $s$ while the leg $t$  becomes the new head (see Figure \ref{fig:varphi-c}). The spanning tree $T$ is unchanged.\vspace{.2cm}\\
\ite For a word  $w=a_1a_2\ldots a_{n}$ on the alphabet $\{a,b,c\}$, we denote by $\varphi_{w}$ the mapping $\varphi_{a_1}\circ\varphi_{a_2}\circ \cdots \circ\varphi_{a_{n}}$.\\
\end{Def}

\begin{Def}\label{def:Phi-prime}
The image of an excursion $w$ by the mapping $\Phi$ is the map with a distinguished spanning tree and a marked external edge obtained by closing the tree-growing map $\varphi_{w}(M_\bullet^0)$, that is, by gluing the head and the unique remaining leg into a marked edge.
\end{Def}



The mapping $\Phi$ has been applied to the excursion $cacbaaccaaba$ in Figure \ref{fig:example-extended-bijection} and \ref{fig:final-step}. Of course, we still need to prove that the mapping $\Phi$ is well defined.

\begin{prop} \label{thm:phi-bien-def}
The mapping $\Phi$ is well defined on any excursion $w$:\\
\ite It is always possible to apply the mapping $\varphi_c$ when required. \\
\ite The tree-growing map $\varphi_{w}(M_\bullet^0)$  has exactly one leg beside the head. This leg and the head are both in the head-face, hence can be glued together. 
\end{prop}

Before proving Proposition \ref{thm:phi-bien-def}, we need three technical results.

\begin{lemma}\label{thm:phi-is-growing}
Let $w$ be a word on the alphabet $\{a,b,c\}$ such that $\varphi_w(M_\bullet^0)$ is well defined. Then, $\varphi_w(M_\bullet^0)$ is a tree-growing map.
\end{lemma}

\dem Let $M_T=\varphi_w(M_\bullet^0)$. It is clear by induction that $T$ is a spanning tree. The only point to prove is that the legs of  $\varphi_w(M_\bullet^0)$ are in the head-face. We proceed by induction on the length of $w$. This property holds for the empty word. If the property holds for $M_T=\varphi_w(M_\bullet^0)$ it clearly holds for $\varphi_a(M_T)$ and $\varphi_b(M_T)$. If $\varphi_c$ can be applied, the head is glued either to the first or to the last leg of $M_T$. Thus, all the remaining legs (including the head of $\varphi_c(M_T)$) are in the same face.
\findem

We shall see shortly (Lemma \ref{thm:all-on-the-T-path}) that whenever the tree-growing map $\varphi_w(M_\bullet^0)$ is well defined, the endpoints of any leg is an ancestor of the head-vertex. Observe that in this case the endpoints of the legs are comparable.

\begin{lemma}\label{thm:deepest-is-first-or-last}
Let $M_T$ be a tree-growing map. Suppose that the endpoint of any leg is an ancestor of the head-vertex. Suppose also that  the first and last legs exist and have distinct endpoints. We call these endpoints $u$ and $v$ with the convention that $u$ is an ancestor of $v$. Then, $v$ is the last vertex incident to a leg on the $T$-path from the root-vertex to the head-vertex.
\end{lemma}

\begin{figure}[h!]
\begin{center}
\input{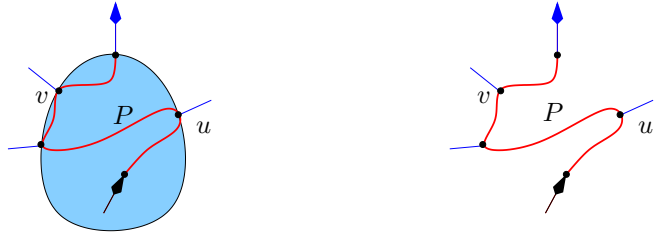}
\caption{The last vertex incident to a leg on the $T$-path from the root-vertex to the head-vertex is $v$.}\label{fig:deepest-is-first-or-last}
\end{center}
\end{figure}

\dem The situation is represented in Figure \ref{fig:deepest-is-first-or-last}. We make an induction on the number of edges that are not in the $T$-path $P$ from the root-vertex to the head-vertex. The property is clearly true if the tree-growing map is reduced to the path $P$ plus some legs. If not, the deletion of a edge not in $P$ does not change the order of appearance of the legs around the head-face. In particular, the first and last legs are unchanged.
\findem

\begin{lemma}\label{thm:all-on-the-T-path}
Let $w$ be a word on the alphabet $\{a,b,c\}$ such  that $\varphi_w(M_\bullet^0)$ is defined. Then the endpoint of any leg of  $\varphi_w(M_\bullet^0)$ is an ancestor of the head-vertex. 
\end{lemma}

\dem We proceed by induction on the length of $w$. The property holds for the empty word. We suppose that it holds for $M_T=\varphi_w(M_\bullet^0)$.  It is clear that the property holds for the tree-growing maps $\varphi_a(M_T)$ and $\varphi_b(M_T)$. If $\varphi_c$ can be applied,  the endpoints of the first and last leg are distinct and comparable. We call these endpoints $u$ and $v$ with the convention that $u$ is an ancestor of $v$. 
By the induction hypothesis, the conditions of Lemma  \ref{thm:deepest-is-first-or-last} are satisfied by $M_T$. Therefore, the vertex $v$ is the last vertex incident to a leg on the $T$-path from the root-vertex to the head-vertex. Hence, any endpoint of a leg of $\varphi_c(M_T)$ is an ancestor of $v$ which is the head-vertex of $\varphi_c(M_T)$.
\findem

\noindent \textbf{Proof of Proposition \ref{thm:phi-bien-def}:}
Let $w$  be an excursion. We consider a suffix  $w'$ of $w$ and denote by $M_T'=\varphi_{w'}(M_\bullet^0)$ the corresponding  tree-growing map (if it is well defined).\vspace{.2cm}\\
 \ite \emph{If $M_T'$ is well defined, it has  $|w'|_a+|w'|_b-2|w'|_c+1$ legs besides the head}. (Observe that, by  \Ref{eq:prefix-extended} and \Ref{eq:word-extended}, the quantity $|w'|_a+|w'|_b-2|w'|_c$ is non-negative.)
\\ 
We proceed by induction on the length of $w'$.  The property holds for the empty word. Moreover, applying $\varphi_a$ or $\varphi_b$ increases by 1 the number of legs whereas applying $\varphi_c$ decreases this number by 2. Thus, the property follows easily by induction.\vspace{.2cm}\\
\ite \emph{The tree-growing map $M_T'$ is well defined}.\\
We proceed by induction on the length of $w'$. The property holds for the empty word. We write $w'=\alpha w''$ and suppose that $M_T''=\varphi_{w''}(M_\bullet^0)$ is well defined.  If $\alpha=a$ or $b$ the tree-growing map $M_T'=\varphi_{\alpha}(M_T'')$ is well defined. We suppose now that  $\alpha=c$. The tree-growing map $M_T''$ has  $|w''|_a+|w''|_b-2|w''|_c+1=|w'|_a+|w'|_b-2|w'|_c+3>2$ legs besides the head. 
 It is clear by induction that all these legs  have distinct endpoints. Moreover, by Lemma \ref{thm:all-on-the-T-path}, all the endpoints of these legs are ancestors of the head-vertex. Thus the endpoints of the legs are comparable. In particular, the endpoints of the first and last legs are comparable. Hence, the mapping  $\varphi_c$ can be applied.\vspace{.2cm}\\
\ite \emph{The tree-growing map $M_T=\varphi_{w}(M_\bullet^0)$ is well defined and has exactly one leg beside the head}.\\
This property follows from the preceding points since $|w|_a+|w|_b-2|w|_c=0$. 
\findem


We now state the key result of this paper. 

\begin{thm}\label{thm:extended-bijection}
The mapping $\Phi$ is a bijection between excursions of size $n$ and  bridgeless 2-near-cubic marked-depth-maps of size~$n$. 
\end{thm}

The proof of Theorem \ref{thm:extended-bijection} is postponed to the next section. For the time being we explore its enumerative consequences. We denote by  $d_n$ the number of bridgeless 2-near-cubic depth-maps of size $n$. Consider a 2-near-cubic map $M$ of size $n$ ($3n+1$ edges, $2n+1$ vertices) and a spanning tree $T$.  Since $T$ has $2n+1$ vertices, $M_T$ has $2n$ internal edges and $n+1$ external edges. Hence, there are  $(n+1)d_n$ bridgeless 2-near-cubic marked-depth-maps. By Theorem \ref{thm:extended-bijection}, this number is equal to the number $e_n$ of excursions of size $n$. Using Proposition \ref{thm:count-extended}, we obtain the following result.

\begin{cor} \label{thm:count-depth-maps}
There are $\displaystyle d_n~=~\frac{e_n}{n+1}~=~\frac{4^n}{(n+1)(2n+1)}{3n \choose n}$ bridgeless 2-near-cubic depth-maps of size $n$. 
\end{cor}

Observe that $d_n$ is also the number of bridgeless cubic depth-maps of size $n-1$ since the bijection between cubic maps and 2-near-cubic maps represented in Figure \ref{fig:2-near-cubic} can be turned into a bijection between cubic depth-maps and 2-near-cubic depth-maps.\\

\section{\textbf{Why the mapping  $\Phi$ is a bijection}}\label{section:preuve-bijection}
In this section, we prove that the mapping $\Phi$ is a bijection between excursions and  bridgeless 2-near-cubic marked-depth-maps. We first prove that the image of any excursion by the mapping $\Phi$ is a bridgeless 2-near-cubic marked-depth-map (Proposition \ref{thm:phi-returns-cubic}). Then we define a mapping $\Psi$ from  bridgeless 2-near-cubic marked-depth-maps to excursions (Definition \ref{def:Psi}) and prove that $\Phi$ and $\Psi$ are inverse mappings (Propositions \ref{thm:inverse1} and \ref{thm:inverse2}).


\begin{prop} \label{thm:phi-returns-cubic}
The image $\Phi(w)$ of any excursion $w$ is a  bridgeless 2-near-cubic marked-depth-map.
\end{prop}

\begin{figure}[h!]
\begin{center}
\input{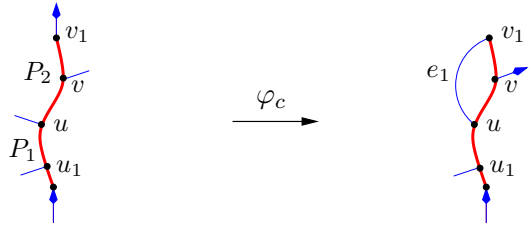}
\caption{Isthmuses are in the $T$-path between $u_0$ and the head-vertex.}\label{fig:watching-isthmuses}
\end{center}
\end{figure}

\dem
Let  $w'$ be a suffix of $w$ and  let $M_T'=\varphi_{w'}(M_\bullet^0)$ be the corresponding  tree-growing map.\vspace{.2cm}\\
\ite \emph{The tree-growing map $M_T'$  is 2-near-cubic}.\\
Applying $\varphi_a$ or $\varphi_b$ creates a new vertex of degree 3 and does not change the degree of the other vertices. Applying $\varphi_c$ does not affect the degree of the vertices. The property follows by induction.\vspace{.2cm}\\
\ite \emph{The head and the root of  $M_T'$ are distinct half-edges.}\\
The property holds for the empty word.  We now write $w'=\alpha w''$.  If $\alpha=a$ or $b$ the property clearly holds for $w'$.  Suppose now that $\alpha=c$.  Let $u$ and $v$ be the vertices incident to the first and last legs of $M_T''=\varphi_{w''}(M_\bullet^0)$ with the convention that $u$ is an ancestor of $v$. By definition, $v$ is the head-vertex of $M_T'=\varphi_c(M_T'')$ and is a proper descendant of $u$. Hence, the head-vertex $v$ and the root-vertex of $M_T'$ are distinct. \vspace{.2cm}\\
\ite \emph{The tree $T$ is a depth tree of  $M_T'$}.\\
The external edges are created by applying the mapping $\varphi_c$, that is, by gluing the head to another leg.  By Lemma \ref{thm:all-on-the-T-path}, any vertex incident to a leg is an ancestor of the head-vertex. Hence, any external edge joins comparable vertices. Moreover, by the preceding point, if the root is part of a complete edge, then this edge is external (internal edges are created by the mappings $\varphi_a$ or $\varphi_b$ which replace the head by a complete edge). \vspace{.2cm}\\
\ite Let $u_0$ be the first vertex of $M_T'$ incident to a leg on the $T$-path from the root-vertex to the head-vertex. \emph{Any isthmus of $M_T'$ is in the $T$-path  between $u_0$ and the head-vertex}.\\
We proceed by induction on the length of $w'$. The property holds for the empty word.  We write $w'=\alpha w''$ and suppose that   it holds for $M_T''=\varphi_{w''}(M_\bullet^0)$.  If $\alpha=a$ or $b$ the property clearly holds for $M_T'=\varphi_{\alpha}(M_T'')$.  We suppose now that $\alpha=c$.  We denote by $u_1$ the first vertex of $M_T''$ incident to a leg on the $T$-path from the root-vertex to the head-vertex. Let $u$ and $v$ be the vertices incident to the first and last legs of $M_T''$ with the convention that $u$ is an ancestor of $v$. By Lemma \ref{thm:all-on-the-T-path}, the vertices $u_1$, $u$ and $v$ are all ancestors of the head-vertex $v_1$ of $M_T''$. Hence, $u$ and $v$ are on the $T$-path between $u_1$ and $v_1$. This situation is represented in Figure \ref{fig:watching-isthmuses}. By definition, the tree-growing map $M_T'$ is obtained from  $M_T''$ by creating an edge $e_1$ between $u$ and $v_1$ while $v$ becomes the new head-vertex.  We denote by $P_1$ (resp. $P_2$) the $T$-path between $u_1$ and $u$ (resp. $u$ and $v_1$). We consider an isthmus $e$ of $M_T'$. The edge $e$ is an isthmus of $M_T''$ (since $M_T''$ is obtained from $M_T$ by deleting an edge). By the induction hypothesis, the isthmus $e$ is either in $P_1$ or in $P_2$. The edge $e$ is not in the path $P_2$ since the new edge $e_1$ creates a cycle with $P_2$. The isthmus  $e$ is in $P_1$, therefore the vertices $u_1$ and $u$ are distinct. Hence $u_1=u_0$ is the first vertex of  $M_T'$ incident to a leg on the $T$-path from the root-vertex to the head-vertex. Thus, the isthmus $e$ is in the $T$-path from $u_0$ to the head-vertex $v$ of $M_T'$.\vspace{.2cm}\\
\ite  \emph{The depth-map $\Phi(w)$ has no isthmus.}\\
By the preceding points, any isthmus of $M_T=\varphi_w(M_\bullet^0)$ is on the $T$-path between the head-vertex and the endpoint of the only remaining leg. Hence, no isthmus remains once the map closed.
\findem

We will now define a mapping $\Psi$ (Definition \ref{def:Psi}) that we shall prove to be the inverse of $\Phi$. The mapping $\Psi$  destructs the tree-growing map that $\Phi$ constructs and recovers the walk. Looking at Figure \ref{fig:example-extended-bijection} from bottom-to-top and right-to-left we see how $\Psi$ works.\\

We first define three mappings $\psi_a$, $\psi_b$, $\psi_c$ on tree-growing maps that we shall prove to be the inverse of $\varphi_a$, $\varphi_b$ and $\varphi_c$ respectively. 
We consider the following conditions for a tree-growing map $M_T$: \\
$(a)$ The head-vertex has degree 3 and is incident to an edge and a leg at the left of the head.\\
$(b)$ The head-vertex has degree 3 and is incident to an edge and a leg at the right of the head.\\
$(c)$ The head-vertex has degree 3 and is incident to 2 edges which are not isthmuses. Furthermore, the tree $T$ is a depth tree.\\

The conditions $(a),~(b),~(c)$ are the domain of definition of  $\psi_a$, $\psi_b$, $\psi_c$ respectively. Before defining these mappings we need a technical lemma.
\begin{lemma}\label{thm:condition-c}
If Condition $(c)$ holds for the tree-growing map $M_T$, then there exists a unique external edge $e_0$ incident to the head-face with one endpoint $u$ ancestor of the head-vertex and one endpoint $v_0$ descendant of the head-vertex. 
\end{lemma}

Lemma \ref{thm:condition-c} is illustrated  by Figure \ref{unique-edge-c}. \\
\begin{figure}[h!]
\begin{center}
\input{unique-edge-c.pstex_t}
\caption{The unique edge $e_0$ satisfying the conditions of Lemma \ref{thm:condition-c}.}\label{unique-edge-c}
\end{center}
\end{figure}

\dem
We suppose that $M_T$ satisfies Condition $(c)$. One of the two edges incident to the head-vertex is in the $T$-path from the root-vertex to the head-vertex. Denote it $e$.  The edge $e$ separates the tree $T$ in two subtrees $T_1$ and $T_2$. We consider the set $E_0$ of external edges having one endpoint in $T_1$ and the other in $T_2$. Any edge satisfying the conditions of Lemma \ref{thm:condition-c} is in $E_0$. Since $e$ is not an isthmus, the set $E_0$ is non-empty. Moreover, any edge in $E_0$ has one endpoint that is a descendant of the head-vertex. Since $T$ is a depth tree, the other endpoint is an ancestor of the head-vertex. It remains to show that there is a unique edge $e_0$ in $E_0$ incident to the head-face. 
By contracting every edge in $T_1$ and $T_2$ we obtain a map with 2 vertices. The edges incident to both vertices are precisely the edges in $E_0\cup \{e\}$. It is clear that exactly 2 of these edges are incident to the head-face.  One is the internal edge $e$ and the other is an external edge $e_0\in E_0$. This edge $e_0$ is the only external edge satisfying the conditions of Lemma \ref{thm:condition-c}.
\findem

We are now ready to define the mappings $\psi_a$, $\psi_b$ and $\psi_c$.
\begin{Def}
Let $M_T$ be a  tree-growing map.\vspace{.2cm} \\
\ite The tree-growing map $M'_{T'}=\psi_a(M_T)$ (resp. $\psi_b(M_T)$) is defined if Condition $(a)$ (resp. $(b)$) holds. 
In this case,  the tree-growing map $M'_{T'}$ is obtained by suppressing the head-vertex $v$ and the 3 incident half-edges. The other half of the edge incident to $v$ becomes the new head.\vspace{.2cm}\\
\ite The tree-growing map $M'_{T'}=\psi_c(M_T)$ is defined if Condition $(c)$  holds.
In this case, we consider the unique external edge $e_0$ with endpoints $u,~v_0$ satisfying the conditions of Lemma \ref{thm:condition-c}. The edge $e_0$ is broken into two legs. The leg incident to $v_0$ becomes the new head (the former head becomes an anonymous leg).\vspace{.2cm}\\
\ite For a word $w=a_1a_2\ldots a_{n}$ on the alphabet $\{a,b,c\}$, we denote by $\psi_{w}$ the mapping $\psi_{a_n}\circ\psi_{a_{n-1}}\circ \cdots \circ\psi_{a_{1}}$.
Moreover, we say that the word $w$ is \emph{readable} on a tree-growing map $M_T$  if the mapping $\psi_{w}$ is well defined on $M_T$.\\
\end{Def}

\noindent \textbf{Remarks:}\\
\ite Applying one of the mappings $\psi_a$, $\psi_b$ or $\psi_c$ to a 2-near-cubic map cannot delete the root (only half-edges incident to a vertex of degree 3 can disappear by application of $\psi_a$ or $\psi_b$).\vspace{.2cm}\\
\ite The conditions $(a)$, $(b)$, $(c)$ are incompatible. Thus, for any tree-growing map $M_T$, there is at most one readable word of a given length. \vspace{.2cm}\\
\ite Applying the mapping $\psi_a$, $\psi_b$ or $\psi_c$  decreases by one the number of edges. Therefore, the length of any  readable word on a tree-growing map $M_T$ is less than or equal to the number of edges in $M_T$.    \\

We now define the mapping $\Psi$ on  bridgeless 2-near-cubic marked-depth-maps. Let $M_T$ be such a map and let $e$ be the marked (external) edge. Observe first that, unless $M_T$ is reduced to a loop, the edge $e$ has two distinct endpoints (or the endpoint of $e$ would be incident to an isthmus). We denote by  $u$ and $v$ the endpoints of $e$ with the convention that $u$ is an ancestor of $v$.  We \emph{open} this map if we disconnect the edge $e$ into two legs and choose the leg incident to $v$ to be the head. We denote by  $M_T^{\dashv \vdash}$ the tree-growing map obtained by opening $M_T$. By convention, opening the  2-near-cubic marked-depth-map reduced to a loop gives $M_\bullet^0$. Note that we obtain $M_T$ by \emph{closing} $M_T^{\dashv \vdash}$. We now define the mapping $\Psi$.

\begin{Def}\label{def:Psi}
Let $M_T$ be a bridgeless 2-near-cubic marked-depth-map. The word $\Psi(M_T)$ is the longest word readable on $M_T^{\dashv \vdash}$.
\end{Def}

We want to prove that $\Phi$ and $\Psi$ are inverse mappings. We begin by proving that the mapping $\psi_\alpha$  is the inverse of $\varphi_\alpha$ for $\alpha=a,b,c$. \\

We say that a tree-growing map satisfies Condition $(c')$ if it satisfies Condition $(c)$ and is such that the endpoint of every leg is an ancestor of the head-vertex.

\begin{lemma}\label{thm:step-by-step-inverse} ~\\
\ite For $\alpha=a$ or $b$,  the mapping $\psi_\alpha \circ \varphi_\alpha$ is the identity on all tree-growing maps and  the mapping $\varphi_\alpha \circ \psi_\alpha$ is the identity on tree-growing maps satisfying Condition $(\alpha)$.\vspace{.2cm}\\
\ite  The mapping $\psi_c \circ \varphi_c$ is the identity on tree-growing maps such that the endpoints of the first and last legs exist and are distinct ancestors of the head-vertex. The mapping $\varphi_c \circ \psi_c$ is the identity on tree-growing maps satisfying Condition $(c')$.
\end{lemma}

Before proving Lemma \ref{thm:step-by-step-inverse}, we need the following technical result.

\begin{lemma}\label{thm:condition-c-prime}
Let $M_T$ be a tree-growing map satisfying Condition $(c')$ and let $e_0$ be the edge with endpoints $u,~v_0$ satisfying the conditions of Lemma  \ref{thm:condition-c}. By definition, the  tree-growing map  $\psi_c(M_T)$ is obtained by breaking $e_0$ into two legs  $s$ and $h$ incident to $u$ and $v_0$ respectively while $h$ becomes the new head. The pair of first and last legs of $\psi_c(M_T)$ is the pair $\{s, t\}$, where $t$ is the head of $M_T$.\end{lemma}

\begin{figure}[h!]
\begin{center}
\input{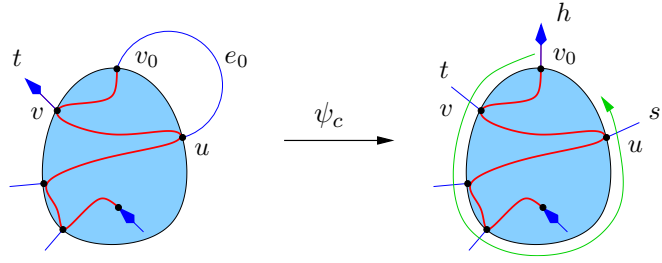}
\caption{The pair of first and last legs of the tree-growing map $\psi_c(M_T)$ is the pair $\{s,t\}$.}\label{fig:destruct-c}
\end{center}
\end{figure}

Lemma \ref{thm:condition-c-prime} is illustrated by  Figure \ref{fig:destruct-c}.\\
 
\noindent \textbf{Proof of Lemma \ref{thm:condition-c-prime}:}\\
\ite Let $v$ be the head-vertex of $M_T$ (i.e. the endpoint of $t$). By Condition $(c')$, the endpoint of any leg of $M_T$ is an ancestor of $v$. Therefore, in the tree-growing map $\psi_c(M_T)$, the vertex $v$ is the last vertex incident to a leg on the $T$-path from the root-vertex to the head-vertex $v_0$. Hence, by Lemma  \ref{thm:deepest-is-first-or-last}, the leg $t$ is either the first or the last leg of $\psi_c(M_T)$.\vspace{.2cm} \\  
\ite No leg lies between $s$ and $h$ on the tour of the head-face of $\psi_c(M_T)$ since this leg would have been inside a non-head face of  $M_T$. Thus the leg $s$ is either the first or the last leg of $\psi_c(M_T)$. 
\findem

\noindent \textbf{Proof of Lemma \ref{thm:step-by-step-inverse}:}\\
\ite For $\alpha=a$ or $b$,  it is clear from the definitions that $\varphi_\alpha \circ \psi_\alpha$ is the identity mapping on all tree-growing maps and that $\varphi_\alpha \circ \psi_\alpha$ is the identity on tree-growing maps satisfying Condition $(\alpha)$.\vspace{.2cm}\\
\ite Consider a tree-growing map $M_T$ such that the endpoints of the first and last legs exist and are distinct ancestors of the head-vertex $v_0$. We call these legs $s$ and $t$ with the convention that the endpoint $u$ of $s$ is an ancestor of the endpoint $v$ of $t$. By definition,  $\varphi_c(M_T)$ is obtained by gluing the head of $M_T$ to $s$ while $t$ becomes the new head. Let $e_0$ be the external edge created by gluing the head to $s$. The head-vertex $v$ of the tree-growing map $\varphi_c(M_T)$ is on the cycle made of $e_0$ and the $T$-path between its two endpoints $u$ and $v_0$, thus  $\varphi_c(M_T)$ satisfies Condition $(c)$. Moreover, the external edge $e_0$ satisfies the conditions of Lemma \ref{thm:condition-c}. Thus, $\psi_c \circ \varphi_c(M_T)=M_T$.\vspace{.2cm}\\
\ite   We consider a tree-growing map $M_T$ satisfying Condition $(c')$.  We consider the edge $e_0$ with endpoints $u,~v_0$ satisfying the conditions of Lemma \ref{thm:condition-c}. By definition,  $\psi_c(M_T)$ is obtained by breaking $e_0$ into two legs  $s$ and $h$ incident to $u$ and $v_0$ respectively while $h$ becomes the new head. By Lemma \ref{thm:condition-c-prime}, the pair of  first and last legs of $\psi_c(M_T)$ is $\{s,t\}$. Moreover, the  endpoint $u$ of $s$ is an ancestor of the endpoint $v$ of $t$ (by definition of $e_0,u,v_0$ in Lemma \ref{thm:condition-c}).  Therefore, the identity $\varphi_c\circ \psi_c(M_T)=M_T$ follows from the definitions.
\findem

\begin{prop}\label{thm:inverse1}
The mapping $\Psi\circ \Phi$ is the identity on excursions.
\end{prop}

\dem\\
\ite \emph{For any word $w$ on the alphabet $\{a,b,c\}$ such  that the tree-growing map $\varphi_w(M_\bullet^0)$ is well defined, the word $w$ is readable on $\varphi_w(M_\bullet^0)$ and $\psi_w \circ \varphi_w(M_\bullet^0)=M_\bullet^0$}.\\
We proceed by induction on the length of $w$. The property holds for the empty word. We write $w=\alpha w'$ with $\alpha=a,b$ or $c$ and suppose that it holds for $w'$. Let $M_T'=\varphi_{w'}(M_\bullet^0)$.  If $\alpha=c$, the  endpoints of the first and last legs of $M_T'$ are distinct and comparable (since $\varphi_c$ is defined on $M_T'$). Moreover, we know by Lemma \ref{thm:all-on-the-T-path} that these endpoints are ancestors of the head-vertex. Thus, for $\alpha=a,b$ or $c$,  Lemma  \ref{thm:step-by-step-inverse} ensures that $\psi_\alpha \circ \varphi_\alpha(M_T')=M_T'$. Therefore, 
$$\psi_{\alpha w'} \circ \varphi_{\alpha w'}(M_\bullet^0)=\psi_{w'}\circ \psi_\alpha \circ \varphi_\alpha \circ \varphi_{w'}(M_\bullet^0)= \psi_{w'}\circ \psi_\alpha \circ \varphi_\alpha (M_T')=\psi_{w'}(M_T'),$$ 
and $\psi_{w'}(M_T')=M_\bullet^0$ by the induction hypothesis. \vspace{.2cm}\\
\ite  \emph{For any excursion $w$, we have $\Psi\circ \Phi(w)=w$}.\\
By definition, the map $M_T=\Phi(w)$ is obtained by closing $\varphi_w(M_\bullet^0)$. 
In order to conclude that $M_T^{\dashv \vdash}=\varphi_w(M_\bullet^0)$, we only need to check that the head of  $M_T^{\dashv \vdash}$ is the head of $\varphi_w(M_\bullet^0)$ (and the non-head leg of $M_T^{\dashv \vdash}$  is the non-head leg of $\varphi_w(M_\bullet^0)$). This is true since the endpoint of the non-head leg of  $\varphi_w(M_\bullet^0)$ is an ancestor of the head-vertex  by Lemma \ref{thm:all-on-the-T-path}. By the preceding point, the word $w$ is readable on $M_T^{\dashv \vdash}=\varphi_w(M_\bullet^0)$ and $\psi_w(M_T)=\psi_w \circ \varphi_w(M_\bullet^0)=M_\bullet^0$. Since no letter is readable on $M_\bullet^0$, the longest word readable on $M_T^{\dashv \vdash}$ is $w$. Thus, $\Psi\circ \Phi(w)=\Psi(M_T)=w$.\\
\findem

It remains to show that $\Phi\circ \Psi$ is the identity mapping on bridgeless 2-near-cubic marked-depth-maps. We first prove that the image of  bridgeless 2-near-cubic marked-depth-maps by $\Psi$ are excursion.

\begin{prop}\label{thm:psi-returns-kreweras}
For any bridgeless 2-near-cubic marked-depth-map $M_T$, the longest word $w$ readable on  $M_T^{\dashv \vdash}$ is an excursion. Moreover, the tree-growing map $\psi_w(M_T^{\dashv \vdash})$ is $M_\bullet^0$.
\end{prop}

\dem If $M_T$ is the map reduced to a loop the result is trivial. We exclude this case in what follows. Let  $w$ be a word readable on $M_T^{\dashv \vdash}$ and  let $N_T=\psi_w(M_T^{\dashv \vdash})$. We denote by $u_0$ the first vertex of $N_T$ incident to a leg on the $T$-path from the root-vertex to the head-vertex.\vspace{.2cm}\\
\ite \emph{Any isthmus of $N_T$ is in the  $T$-path between $u_0$ and the head-vertex}.\\
We proceed by induction on the length of $w$. Suppose first that $w$ is the empty word. Let $e_0$ be the marked edge of $M_T$. By definition, the tree-growing map $N_T=M_T^{\dashv \vdash}$ is obtained from  $M_T$ by breaking $e_0$ into two legs: the head and another leg incident to $u_0$. Let $e$ be an isthmus of $N_T$ and let $N_1,~N_2$ be the two connected submaps obtained by deleting $e$. Since $e$ is not an isthmus of $M_T$, the edge $e_0$ joins $N_1$ and $N_2$. Therefore, the root-vertex and head-vertex are not in the same submap. Thus, the isthmus $e$ is in any path between $u_0$ and the head-vertex, in particular  it is in the $T$-path.\\
We now write $w=\alpha w'$ with $\alpha=a,b$ or $c$ and suppose, by the induction hypothesis, that the property holds for $w'$. We denote by $u_0'$ the first vertex of $N_T'=\psi_{w'}(M_T^{\dashv \vdash})$ incident to a leg on the $T$-path from the root-vertex to the head-vertex. Suppose first that $\alpha=a$ or $b$. The edge incident to the head-vertex of  $N_T'$ is an isthmus hence, by the induction hypothesis, it is in the $T$-path between $u_0'$ and the head-vertex $v_0'$ of $N_T'$.  Hence, $u_0'\neq v_0'$. Thus,  $u_0=u_0'$ and every isthmus  of $N_T=\psi_\alpha(N_T')$ is in the  $T$-path between $u_0$ and the head-vertex. 
Suppose now that $\alpha=c$. Since $w$ is readable on $M_T^{\dashv \vdash}$, the tree-growing map  $N_T'=\psi_{w'}(M_T^{\dashv \vdash})$ satisfies Condition $(c)$.  We consider the edge $e_0$ with endpoints $u,~v_0$ satisfying the conditions of Lemma \ref{thm:condition-c}. The map $N_T=\psi_{c}(N_T')$ is obtained from $N_T$ by breaking $e_0$ into two legs. By definition, the head-vertex of $N_T$ is $v_0$. Moreover, the vertex $u_0$ is either $u_0'$ or $u$ if $u$ is an ancestor of $u_0'$. We consider an isthmus $e$ of $N_T$. If $e$ is an isthmus of  $N_T'$, it is in the $T$-path between $u_0'$ to the head-vertex of $N_T'$ which is included in the $T$-path between $u_0$ and  $v_0$. If $e$ is not an isthmus of  $N_T'$, we consider the two connected submaps $N_1,~N_2$ obtained from $N_T$ by deleting the isthmus $e$.  Since $e$ is not an isthmus of $N_T'$, the edge $e_0$ joins $N_1$ and $N_2$. Hence, the endpoints $u$ and $v_0$ of $e_0$ are not in the same submap. Thus, the isthmus $e$ is in every path of $N_T$ between $u$ and the head-vertex $v_0$, in particular, it is in the $T$-path between $u_0$ and $v_0$.\vspace{.2cm} \\
\ite \emph{The tree-growing map $N_T$ has at least one leg beside the head}.\\
We proceed by induction. The property holds for the empty word. We now write $w=\alpha w'$ with $\alpha=a,b$ or $c$ and suppose that the property holds for $w'$. Suppose first that $\alpha=a$ or $b$. Since Condition  $(\alpha)$ holds, the edge incident to the head-vertex $v_0'$ of the tree-growing map  $N_T'=\psi_{w'}(M_T^{\dashv \vdash})$ is an isthmus. By the preceding point, this edge is  on the $T$-path between $u_0'$ and  $v_0'$, where $u_0'$ be the first vertex of $N_T'$ incident to a leg on the $T$-path from the root-vertex to the head-vertex. Thus  $u_0'\neq v_0'$ and $N_T=\psi_\alpha(N_T')$ has at least one leg (the one incident to $u_0'$) beside the head. In the case $\alpha=c$,  the tree-growing map $N_T=\psi_c(N_T')$ has one more legs than $N_T'$, hence it has at least one leg beside the head.\vspace{.2cm}\\
\ite \emph{The head and root of $N_T$ are distinct half-edges}. \\
By definition, the map $M_T^{\dashv \vdash}$ has one leg beside the head whose endpoint is a proper ancestor of the head-vertex. Hence, the head-vertex and root-vertex are distinct. We suppose now that  $w=\alpha w'$ with $\alpha=a,b$ or $c$. If $\alpha=a$ or $b$ the head of $N_T$ is an half-edge of  $N_T'=\psi_{w'}(M_T^{\dashv \vdash})$ which is part of an internal edge. Hence it is not the root. If $\alpha=c$, the head of $N_T$  is part of an external edge $e$ of  $N_T'=\psi_{w'}(M_T^{\dashv \vdash})$. The edge is broken into the head of $N_T$ and another leg whose endpoint is a proper ancestor of the head-vertex. Hence, the head-vertex and root-vertex of $N_T$ are distinct.  \vspace{.2cm}\\
\ite \emph{If $w$ is the longest readable word, then $N_T=M_\bullet^0$}.\\
We first prove that the root-vertex and the head-vertex of $N_T$ are the same. Suppose they are distinct. In this case, the head-vertex has degree 3 and is incident to at least one edge. If it is incident to one edge, then one of the conditions $(a)$ or $(b)$ holds and $w$ is not the longest readable word.  Hence the head-vertex is incident to two edges $e_1$ and $e_2$.  One of these edges, say $e_1$, is in the $T$-path from the root-vertex to the head-vertex and the other $e_2$ is not. By a preceding point, the edge $e_2$ is not an isthmus . Therefore, $e_1$ is not an isthmus either ($e_1$ and $e_2$ have the same ability to disconnect the map). In this case, Condition $(c)$ holds (since $T$ is a depth tree) and  $w$ is not the longest readable word. Thus, the root-vertex and the head-vertex of $N_T$ are the same.  Therefore, the root-vertex  has degree 2 and is incident to the head and the root. The head and the root are distinct (by the preceding point). Moreover the root is a leg. Indeed, if the root was not a leg it would be part of an external edge which is an isthmus (which is impossible since the tree $T$ is spanning). Hence the root-vertex is incident to two legs: the root and the head. Thus, $N_T=M_\bullet^0$. \vspace{.2cm}\\
\ite \emph{The tree-growing map $N_T$ has $2|w|_c-|w|_a-|w|_b+1$ legs beside the head}.\\
The tree-growing map $M_T^{\dashv \vdash}$ has one leg beside the head. Moreover, applying mapping $\psi_a$ or $\psi_b$ decreases by one the number of legs whereas applying mapping $\psi_c$ increases this number by two. Hence the property follows easily by induction. \vspace{.2cm}\\
 \ite  \emph{The longest word $w$ readable on $M_T^{\dashv \vdash}$ is an excursion}.\\
By the preceding points, any prefix $w'$ of $w$ satisfies  $2|w'|_c-|w'|_a-|w'|_b+1\geq 1$ (since this quantity is the number of non-head legs of $\psi_{w'}(M_T^{\dashv \vdash})$). Moreover, since  $\psi_{w}(M_T^{\dashv \vdash})=M_\bullet^0$ has one leg beside the root, we have  $2|w|_c-|w|_a-|w|_b+1=1$.  These properties are equivalent to  \Ref{eq:prefix-extended} and \Ref{eq:word-extended}, hence $w$ is an excursion.
\findem

\begin{prop}\label{thm:inverse2}
The mapping $\Phi\circ \Psi$ is the identity on  bridgeless 2-near-cubic marked-depth-maps. 
\end{prop}

\dem Let $M_T$ be  a  bridgeless 2-near-cubic marked-depth-map.\vspace{.2cm}\\
\ite \emph{For any word $w$ readable on   $M_T^{\dashv \vdash}$, the endpoints of any leg of  $\psi_w(M_T^{\dashv \vdash})$ is an ancestor of the head-vertex.}\\
We proceed by induction on the length of $w$. The property holds for the empty word. We now write $w=\alpha w'$ with $\alpha=a,b$ or $c$ and suppose that it holds for $w'$. For $\alpha=a$ or $b$, the property clearly holds for $w$. Suppose now that $\alpha=c$. Since $w$ is readable, the tree-growing map  $N_T'=\psi_{w'}(M_T^{\dashv \vdash})$ satisfies Condition $(c)$.  We consider the edge $e_0$ with endpoints $u,~v_0$ satisfying the conditions of Lemma \ref{thm:condition-c}. By definition, the head-vertex $v_0$ of $N_T=\psi_{c}(N_T')$ is a descendant of the head-vertex $v$ of $N_T'$. By the induction hypothesis, the endpoint of any leg of $N_T'$ is an ancestor of $v$. Hence, the endpoint of any leg of $N_T$ is an ancestor of the head-vertex $v_0$.\vspace{.2cm}\\
\ite  \emph{For any word $w$ readable on   $M_T^{\dashv \vdash}$, we have $\varphi_w \circ \psi_w(M_T^{\dashv \vdash})=M_T^{\dashv \vdash}$}.\\
We proceed by induction. The property holds for the empty word. We now write $w=\alpha w'$ with $\alpha=a,b$ or $c$ and suppose that the property  holds for $w'$. If $\alpha=a$ or $b$  the induction step is given directly by Lemma \ref{thm:step-by-step-inverse} (since Condition $(\alpha)$ holds for $M_T'=\psi_{w'}(M_T^{\dashv \vdash})$).  If $\alpha=c$, that is, Condition $(c)$ holds for $M_T'=\psi_{w'}(M_T^{\dashv \vdash})$, we must prove that Condition  $(c')$ holds (in order to apply Lemma \ref{thm:step-by-step-inverse}). But we are ensured that Condition $(c')$ holds by the preceding point.
 Thus, for $\alpha=a,b$ or $c$,  Lemma  \ref{thm:step-by-step-inverse} ensures that $\varphi_\alpha \circ \psi_\alpha(M_T')=M_T'$. Therefore, 
$$\varphi_{\alpha w'} \circ \psi_{\alpha w'}(M_T^{\dashv \vdash})=\varphi_{w'}\circ \varphi_\alpha \circ \psi_\alpha \circ \psi_{w'}(M_T^{\dashv \vdash})= \varphi_{w'}\circ \varphi_\alpha \circ \psi_\alpha (M_T')=\varphi_{w'}(M_T'),$$
 and $\varphi_{w'}(M_T')=M_T^{\dashv \vdash}$ by the induction hypothesis. \vspace{.2cm}\\
\ite \emph{$\Phi\circ \Psi(M_T)=M_T$.}\\
By definition, the word $w=\Psi(M_T)$ is the longest readable word on $M_T^{\dashv \vdash}$. Hence, by Proposition \ref{thm:psi-returns-kreweras}, $\psi_w(M_T^{\dashv \vdash})=M_\bullet^0$. By the preceding point, $\varphi_w(M_\bullet^0)=\varphi_w\circ \psi_w(M_T^{\dashv \vdash})=M_T^{\dashv \vdash}$. By definition, the map $\Phi(w)$ is obtained by closing $\varphi_w(M_\bullet^0)=M_T^{\dashv \vdash}$, hence $\Phi(w)=M_T$. Thus, $\Phi\circ \Psi(M_T)=\Phi(w)=M_T$. 
\findem

By Proposition \ref{thm:phi-returns-cubic}, the mapping $\Phi$ associates a  bridgeless 2-near-cubic marked-depth-map with any excursion. Conversely, by Proposition \ref{thm:psi-returns-kreweras}, the mapping $\Psi$ associates an excursion with any  bridgeless 2-near-cubic marked-depth-map. The mappings  $\Phi$ and $\Psi$ are inverse mappings by Propositions \ref{thm:inverse1} and \ref{thm:inverse2}. Thus, the mapping $\Phi$ is a bijection between excursions and  bridgeless 2-near-cubic marked-depth-maps. Moreover, if an excursion $w$ has size $n$ (length $3n$), the  2-near-cubic depth-map $\Phi(w)$ has size $n$ ($3n+1$ edges). This concludes the proof of Theorem~\ref{thm:bijection}.
\findem

\section{\textbf{A bijection between Kreweras walks and cubic depth-maps}} \label{section:bijection}
In this section, we prove that the mapping $\Phi$ establishes a bijection between Kreweras walks ending at the origin and 2-near-cubic depth-maps. This result is stated more precisely in the following theorem.

\begin{thm}\label{thm:bijection}
Let $w$ be an excursion. The marked edge of the 2-near-cubic depth-map $\Phi(w)$ is the root-edge if and only if the excursion $w$ is a  Kreweras walk ending at the origin. \\
Thus, the mapping $\Phi$ induces a bijection between Kreweras walks of size $n$ (length $3n$) ending at the origin and bridgeless 2-near-cubic depth-maps of size $n$ ($3n+1$ edges).\\
\end{thm}

\begin{figure}[h!]
\begin{center}
\input{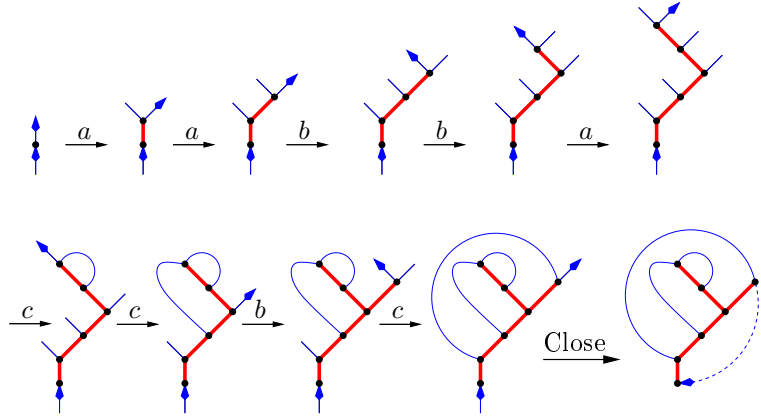}
\caption{The image of a Kreweras walk by $\Phi$: the root-edge is marked.}\label{fig:example-bijection}
\end{center}
\end{figure}

Figure \ref{fig:example-bijection} illustrates an instance of Theorem \ref{thm:bijection}. Before proving this theorem we explore its enumerative consequences. From Theorem \ref{thm:bijection}, the number $k_n$ of  Kreweras walks of size $n$ is equal to the number $d_n$ of bridgeless 2-near-cubic depth-maps of size $n$. The number $d_n$ is given by Corollary \ref{thm:count-depth-maps}. We obtain the following result.

\begin{thm}\label{thm:count-kreweras}
There are $\displaystyle k_n ~=~ \frac{4^n}{(n+1)(2n+1)} {3n \choose n}$ Kreweras walks of size $n$ (length $3n$) ending at the origin. 
\end{thm}

The rest of this section is devoted to the proof of  Theorem \ref{thm:bijection}. \\

Consider a growing map $M$ such that the root is a leg. Recall that \emph{making the tour of the head-face} means following  its border in counterclockwise direction  starting from the head (see Figure \ref{fig:making-the-tour}). We call \emph{left} (resp. \emph{right}) the legs encountered before (resp. after) the root during the tour of the head-face. For instance, the growing map in Figure \ref{fig:making-the-tour} has one left leg and two right legs.\\

\begin{lemma}\label{thm:kreweras-gives-root}
For any Kreweras walk $w$ ending at the origin, the marked edge of  $\Phi(w)$ is the root-edge.
\end{lemma}

\dem 
Let $w'$ be a suffix of $w$ and let $M_T'=\varphi_{w'}(M_\bullet^0)$ be the corresponding tree-growing map. \vspace{.2cm}\\ 
\ite  \emph{The root of $M_T'$  is a leg and  $M_T'$  has   $|w'|_a-|w'|_c$ left legs and  $|w'|_b-|w'|_c$ right legs}. (Observe that, these quantities are non-negative by \Ref{eq:prefix} and \Ref{eq:word}.)\\
We proceed by induction  on the length of $w'$. The property holds for the empty word. We now write $w'=\alpha w''$ with $\alpha=a,b$ or $c$ and suppose that the property  holds for $w''$. If $\alpha=a$ or $b$ the property holds for $w'$ since applying $\varphi_a$ (resp. $\varphi_b$) increases by one the number of left (resp. right) legs. We now suppose that $\alpha=c$. We know that   $|w''|_a-|w''|_c=|w'|_a-|w'|_c+1\geq 1$. Hence, by the induction hypothesis, the tree-growing-map  $M_T''=\varphi_{w''}(M_\bullet^0)$ has at least one left leg. Similarly,  $M_T''$ has at least one right leg.  Therefore, the first (resp. last) leg of $M_T''$ is a left (resp. right) leg. Hence, applying $\varphi_c$ to $M_T''$ decreases by one the number of left (resp. right) legs. Thus, the property holds for $w'$.\vspace{.2cm}\\
\ite For $w'=w$, the preceding point shows that $\varphi_w(M_\bullet^0)$ has only one leg beside the head and that this leg is the root. Thus, the marked edge of $\Phi(w)$ is the root-edge.
\findem

\begin{lemma}\label{thm:root-gives-kreweras}
For any bridgeless 2-near-cubic depth-map $M_T$ marked on the root-edge, the word $w=\Psi(M_T)=\Phi^{-1}(M_T)$ is  a Kreweras walk ending at the origin.
\end{lemma}

\dem
Let  $w$ be a word readable on $M_T^{\dashv \vdash}$ and  let $N_T=\psi_w(M_T^{\dashv \vdash})$. Observe that the root of $N_T$ is a leg (since it is the case in $M_T^{\dashv \vdash}$ and the root never disappears).\vspace{.2cm}\\
\ite \emph{The tree-growing map $N_T$ has $|w|_c-|w|_a$ left legs and  $|w|_c-|w|_b$ right legs}.\\
We proceed by induction  on the length of $w$. The property holds for the empty word.  We now write $w=\alpha w'$ with $\alpha=a,b$ or $c$ and suppose that the property  holds for $w'$. If $\alpha=a$ or $b$ the property holds for $w$ since applying $\psi_a$ (resp. $\psi_b$) decreases by one the number of left (resp. right) legs.  We now suppose that $\alpha=c$. The map $N_T'=\psi_{w'}(M_T^{\dashv \vdash})$ satisfies Condition $(c)$. We have already proved (see the first point in the proof of Lemma \ref{thm:psi-returns-kreweras}) that the endpoint of every leg is an ancestor of the head-vertex. Hence $N_T'$ satisfies Condition $(c')$. Therefore, Lemma \ref{thm:condition-c-prime} holds for $N_T'$. We adopt the notations $h,s,t$ of this lemma which is illustrated in  Figure~\ref{fig:destruct-c}.  By Lemma \ref{thm:condition-c-prime}, the pair of first and last head of $N_T=\psi_c(N_T')$ is the pair $\{s,t\}$. Hence, in the pair $\{s,t\}$ one is a left leg and the other is a right leg of $N_T$. Moreover, the other  left and right legs of  $N_T$ are the same as in $N_T'$.  Thus, applying $\psi_c$ to $N_T'$ increases by one the number of left (resp. right) legs. Hence, the property holds for $w$.\vspace{.2cm}\\
\ite \emph{The word $w=\Psi(M_T)$ is  a Kreweras walk ending at the origin.}\\
By definition,  $w$ is the longest word readable on $M_T^{\dashv \vdash}$. By Proposition \ref{thm:psi-returns-kreweras}, $\psi_{w}(M_T^{\dashv \vdash})=M_\bullet^0$. By the preceding point, we get   $|w|_c-|w|_a=0$ and  $|w|_c-|w|_b=0$ (since $M_\bullet^0$ has no left nor right leg). Moreover, for any suffix $w'$ of $w$, the preceding point proves that  $|w'|_c-|w'|_a\geq 0$ and  $|w'|_c-|w'|_b\geq 0$. These properties are equivalent to \Ref{eq:prefix} and \Ref{eq:word}, hence $w$ is  a Kreweras walk ending at the origin.
\findem

\section{\textbf{Enumerating depth trees and cubic maps}} \label{section:counting-cubic}
In Section \ref{section:extension-bijection}, we exhibited a bijection $\Phi$ between excursions and  bridgeless 2-near-cubic marked-depth-maps. As a corollary we obtained the number of bridgeless 2-near-cubic depth-maps of size $n$:  
$d_n~=~ \frac{4^n}{(n+1)(2n+1)} {3n \choose n}.$
In this section, we prove that any bridgeless 2-near-cubic map of size $n$ has $2^n$ depth trees (Corollary \ref{thm:counting-depth-tree-maps}). Hence, the number of bridgeless 2-near-cubic maps of size $n$ is $ c_n=\frac{d_n}{2^n}=\frac{2^n}{(n+1)(2n+1)} {3n \choose n}.$ 
Given the bijection between 2-near-cubic maps and cubic maps (see Figure \ref{fig:2-near-cubic}), we obtain the following theorem.
\begin{thm}\label{thm:number-cubic}
There are $\displaystyle c_n ~=~ \frac{2^n}{(n+1)(2n+1)} {3n \choose n}$ bridgeless cubic maps with $3n$ edges.  
\end{thm}

By duality, $c_n$ is also the number of loopless triangulations with $3n$ edges. Hence, we recover Equation \Ref{eq:triangulations} announced in the introduction. As mentioned above, an alternative bijective proof of Theorem \ref{thm:number-cubic} was given in~\cite{Schaeffer:triangulation}. \\

The rest of this section is devoted to the counting of depth trees on cubic maps and, more generally, on cubic (potentially non-planar) graphs. We first give an alternative characterization of depth trees. This characterization is based on the depth-first search (DFS) algorithm (see Section 23.3  of \cite{Cormen:introduction-algorithms}). We consider the DFS algorithm as an algorithm for constructing a spanning tree of a graph. \\


Consider a graph $G$ with a distinguished vertex $v_0$. If the DFS algorithm starts at $v_0$, the subgraph $T$ (see below) constructed by the algorithm remains a tree containing $v_0$.  We call \emph{visited} the vertices in $T$ and \emph{unvisited} the other vertices. The distinguished vertex $v_0$ is considered as the root-vertex of the tree. Hence, any vertex in $T$ distinct from $v_0$ has a \emph{father} in $T$. 

\begin{Def} \label{def:DFS}
Depth-first search (DFS) algorithm.  \vspace{.3cm}\\
\textbf{Initialization:} Set the \emph{current vertex} to be $v_0$ and the tree $T$ to be reduced to $v_0$.\vspace{.1cm}\\
\textbf{Core:} While the current vertex $v$ is adjacent to some unvisited vertices or is distinct from $v_0$ do:\\
If there are some edges linking the current vertex $v$ to an unvisited vertex, then choose one of them. Add the chosen edge $e$ and its unvisited endpoint $v'$ to the tree $T$. Set  the current vertex to be $v'$.\\
Else, backtrack, that is, set the current vertex to be the father of $v$ in~$T$. \vspace{.1cm}\\
\textbf{End:} Return the tree $T$.\\
\end{Def}

It is well known that the DFS algorithm returns a spanning tree. It is also known \cite{Cormen:introduction-algorithms} that the two following properties are equivalent for a spanning tree $T$ of a graph $G$ having a distinguished vertex $v_0$:\\
$\textbf{(i)}$ Any external edge joins comparable vertices. \\
$\textbf{(ii)}$ The tree $T$ can be obtained by a DFS algorithm on the graph $G$ starting at $v_0$.\\

Before stating the main result of this section, we need an easy preliminary lemma. 

\begin{lemma}
Let $G$ be a connected graph with a distinguished vertex $v_0$ whose deletion does not disconnect the graph. Then, any spanning tree $T$ of $G$ satisfying conditions \textbf{(i)}-\textbf{(ii)} has exactly one edge incident to $v_0$. \\
\end{lemma}

\dem
Let $e_0$ be an edge of $T$ incident to $v_0$ and let $v_1$ be the other endpoint of $e_0$.  We partition the vertex set $V$ of $G$ into $\{v_0\}\cup V_0\cup V_1$, where  $V_1$ is the set of descendants of $v_1$.   There is no internal edge joining a vertex in $V_0$ and a vertex in $V_1$. There is no external edge either or it would join two non-comparable vertices. Thus $V_0=\emptyset$ or the deletion of $v_0$ would disconnect the graph.
\findem

\begin{thm} \label{thm:counting-depth-tree}
Let $G$ be a loopless connected graph with a distinguished vertex $v_0$ whose deletion does not disconnect the graph. 
Let $e_0$ be an edge incident to $v_0$. If $G$ is a $k$-near-cubic  graph ($v_0$ has degree $k$ and the other vertices have degree $3$) of size $n$ ($3n+2k-3$ edges), then there are $2^n$ trees containing $e_0$ and satisfying conditions \textbf{(i)}-\textbf{(ii)}.
\end{thm}

Given that the depth trees are the spanning trees satisfying conditions \textbf{(i)}-\textbf{(ii)} and not containing the root, the following corollary is immediate.

\begin{cor} \label{thm:counting-depth-tree-maps}
Any bridgeless 2-near-cubic map of size $n$ ($3n+1$ edges) has $2^n$ depth trees.
\end{cor}

\noindent \textbf{Remark:} Theorem \ref{thm:counting-depth-tree} implies that any $k$-near-cubic loopless graph of size $n$ has $k2^n$  trees satisfying the conditions \textbf{(i)}-\textbf{(ii)}.\\

The rest of this section is devoted to the proof of Theorem \ref{thm:counting-depth-tree}. The proof relies on the intuition that exactly $n$ \emph{real} binary choices have to be made during the execution of a DFS algorithm on a $k$-near-cubic map of size $n$.\\

Given a graph $G$ and a subset of vertices $U$, we say that two vertices $u$ and $v$ are \emph{$U$-connected} if  there is a path between $u$ and $v$ containing only vertices in $U\cup \{u,v\}$. 

\begin{lemma} \label{thm:visited-vertices}
Let $v$ be the current vertex and let $U$ be the set of unvisited vertices at a given time of the DFS algorithm. The vertices that will be visited before the last visit to $v$ are the vertices in $U$ that are $U$-connected to $v$. 
\end{lemma}

\begin{figure}[h!]
\begin{center}
\input{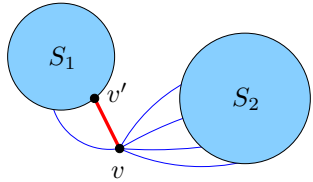}
\caption{Partition of the vertices in $S$.}\label{fig:unvisited-vertices}
\end{center}
\end{figure}

\dem
Let $S$ be the set of vertices in  $U$ that are $U$-connected to $v$. We make an induction on the cardinality of $S$.  If the set $S$ is empty, there is no edge linking $v$ to an unvisited vertex. Hence, the next step in the algorithm is to backtrack and the vertex $v$ will never be visited again. In other words, it is the last visit to $v$, hence the property holds. Suppose now that $S$ is non-empty. In this case, there are some edges linking the current vertex $v$ to an unvisited vertex. Let $e$ be the edge chosen by the DFS algorithm and let $v'\in U$ be the corresponding endpoint. Let  $S_1$ be the set of vertices in  $U$ that are $U$-connected to $v'$ and let $S_2=S\setminus S_1$. Observe that no edge joins a vertex in $S_1$ and a vertex in $S_2$. This situation is represented in Figure \ref{fig:unvisited-vertices}. The set of vertices in  $U'=U\setminus \{v'\}$ that are $U'$-connected to $v$ is $S_1'=S_1\setminus \{v'\}$ (since a vertex is $U$-connected to $v'$ if and only if it is $U'$-connected to $v'$). By the induction hypothesis, $S_1'$ is the set of vertices visited between the first and last visit to $v'$. Hence $S_1$ is the set of vertices visited before the algorithm returns to $v$. Since no edge joins a vertex in $S_1$ and a vertex in $S_2$, the vertices in $S_2$ are the vertices in $U\setminus S_1$ that are $(U\setminus S_1)$-connected to $v$. By the induction hypothesis, $S_2$ is the set of vertices visited before the last visit to $v$. Thus, the property holds.
\findem

\begin{figure}[h!]
\begin{center}
\input{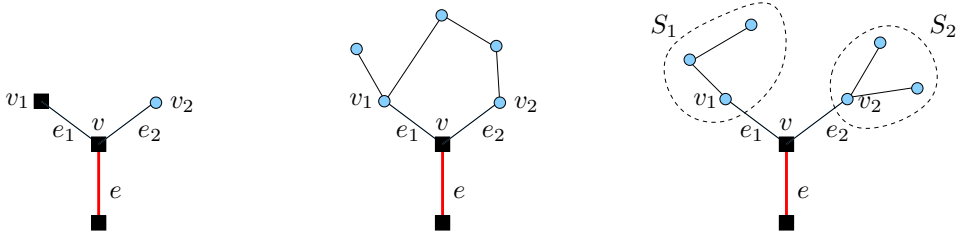}
\caption{Case $(\alpha)$ (left), case $(\beta)$ (middle) and case $(\gamma)$ (right). The visited vertices are indicated by a square while unvisited ones are indicated by a circle.}\label{fig:three-cases}
\end{center}
\end{figure}

\noindent \textbf{Proof of Theorem \ref{thm:counting-depth-tree}:}  
Clearly, the spanning trees containing $e_0$ and satisfying the conditions  \textbf{(i)}-\textbf{(ii)} are the spanning trees obtained by a DFS algorithm for which the first \emph{core step} is to choose $e_0$. We want to prove that there are $2^n$ such spanning trees.\\
We consider an execution of the DFS algorithm  for which the first \emph{core step} is to choose $e_0$ and denote by $\mathcal{T}$ the spanning tree returned by the DFS algorithm (in order to distinguish it from the evolving tree $T$). After the first \emph{core step}, the tree $T$ is reduced to $e_0$ and its two endpoints $v_0$ and $v_0'$. Let $V$ be the vertex set of $G$ and let $V'=V\setminus \{v_0,v_0'\}$. Since the deletion of $v_0$ does not disconnect the graph, every vertex in $V'$ is $V'$-connected to $v_0'$. Hence, by Lemma  \ref{thm:visited-vertices}, every vertex  will be visited before the algorithm returns to $v_0$. Thus, from this stage on, the current vertex $v$ is incident to 3 edges  $e,~ e_1,~ e_2$, where $e\in T$ links $v$ to its father.\vspace{.2cm}\\
\ite We denote by $v_1$ and $v_2$ the  endpoints of $e_1$ and $e_2$ respectively (these endpoints are not necessarily distinct) and we denote by $U$ the set of unvisited vertices. We distinguish three cases:\\
$(\alpha)$ at least one of the vertices $v_1,~v_2$ is not in $U$,\\
$(\beta)$ the two vertices  $v_1,~v_2$ are in $U$ and are  $U$-connected with each other,\\
$(\gamma)$ the two vertices  $v_1,~v_2$  are in $U$ and are not $U$-connected with each other.\\
The three cases are illustrated by Figure \ref{fig:three-cases}. We prove successively the following properties:\\
\iten  \emph{In case $(\alpha)$, no choice has to be done by the algorithm}. \\
Indeed, there is at most one edge ($e_1$ or $e_2$) linking the current vertex $v$ to an unvisited vertex.\\
\iten  \emph{In case $(\beta)$, the algorithm has to choose between $e_1$ and $e_2$. This choice \emph{necessarily} leads to two different spanning trees $\mathcal{T}$. Indeed the edge $e_1$ (resp. $e_2$) is in $\mathcal{T}$ if and only if the choice of $e_1$ (resp. $e_2$) is made.}\\
Suppose (without loss of generality), that the choice of $e_1$ is made. The vertex $v_2$ is $(U\cup \{v_1\})$-connected to $v_1$ (a vertex is $(U\cup \{v_1\})$-connected to $v_1$ if and only if it is $U$-connected to $v_1$). Hence, by Lemma \ref{thm:visited-vertices}, the vertex $v_2$ will be visited before the last visit to $v_1$, that is, before the algorithm returns to the vertex $v$. Therefore, the edge $e_2$ will not be in the spanning tree $\mathcal{T}$.\\
\iten  \emph{In case $(\gamma)$, the algorithm has to choose between $e_1$ and $e_2$.  Moreover, any tree $\mathcal{T}$ obtained by choosing $e_1$ can be also obtained by choosing $e_2$}.\\
Let $S_1$ and $S_2$ be the set of vertices in $U$ that are $U$-connected to $v_1$ and $v_2$ respectively. Observe that the sets $S_1$ and $S_2$ are disjoint and no edge links a vertex in $S_1$ and a vertex in $S_2$ (otherwise the vertices $v_1$ and $v_2$ would be  $U$-connected). Suppose that the choice of $e_1$ is made. The set of vertices in $U\setminus \{v_1\}$ that are $U\setminus \{v_1\}$-connected to $v_1$ is $S_1\setminus \{v_1\}$. Hence, by Lemma \ref{thm:visited-vertices}, the set of vertices visited before the last visit to $v_1$, that is, before the algorithm returns to $v$ is $S_1$. Since $v_2$ is not in $S_1$ the next step of the algorithm is to choose $e_2$. Let $U_2=U\setminus S_1$ be the set of unvisited vertices at this stage. Since  no vertex in $S_1$ is adjacent to a vertex in $S_2$, the set of vertices in $U_2\setminus \{v_2\}$ that are $(U_2\setminus \{v_2\})$-connected to $v_2$ is $S_2\setminus \{v_2\}$. Hence, by Lemma \ref{thm:visited-vertices}, the set of vertices visited before the last visit to $v_2$, that is, before the algorithm returns to $v$ is $S_2$. Let $T_1$ (resp. $T_2$) be the subtree constructed by the algorithm between the first and last visit to $v_1$ (resp. $v_2$). Since no vertex in $S_1$ is adjacent to a vertex in $S_2$, the subtree $T_1$ could have been constructed exactly the same way if the algorithm had chosen $e_2$ (instead of $e_1$) at the beginning. Similarly, the subtree $T_2$ could have been constructed exactly in the same way if the algorithm had chosen $e_2$ at the beginning.  Therefore, the tree  $\mathcal{T}$ returned by the algorithm could have been constructed if the algorithm had chosen $e_2$ (instead of $e_1$) at the beginning.\vspace{.2cm}\\
\ite \emph{During any execution of the DFS algorithm we are exactly $n$ times in case $(\beta)$}. \\
The $k$-near-cubic graph $G$ has $3n+2k-3$ edges and $2n+2k-1$ vertices. Hence, the spanning tree $\mathcal{T}$ has $2n+2k-2$ edges. Thus, there are $n+k-1$ external edges among which $k-1$ are incident to $v_0$. Let $E_\beta$ be the set of the $n$ external edges not incident to $v_0$. Since $G$ is loopless and the spanning tree  $\mathcal{T}$ satisfies  \textbf{(i)}-\textbf{(ii)}, the edges in $E_\beta$ have distinct and comparable endpoints.  For any edge $e$ in $E_\beta$, we denote by $v_e$ the endpoint of $e$ which is the ancestor of the other endpoint. The vertex $v_e$ is incident to $e$, to the edge of $\mathcal{T}$ linking $v$ to its father and to another edge in $\mathcal{T}$ linking $v_e$ to its son (otherwise $v_e$ has no descendant). In particular, if $e$ and $e'$ are distinct edges in $E_\beta$, then the vertices  $v_e$ and $v_{e'}$ are distinct. Thus, the set of vertices $V_\beta=\{v_e / e\in E_\beta \}$ has size $n$.\\
We want to prove that the case $(\beta)$ occurs when the algorithm visit a vertex in $V_\beta$ for the first time (and not otherwise). Let $v$ be a vertex in $V_\beta$. The vertex $v$ is incident to an edge $e_1$ in $E_\beta$, an edge $e$ in $\mathcal{T}$ linking $v$ to its father and another edge $e_2$ in $\mathcal{T}$ linkink $v$ to its son. Let $T$ be the tree constructed by the algorithm at the time of the first visit to $v$ and let $U$ be the set of unvisited vertices. Any descendant of $v$ is in $U$. In particular, the endpoints $v_1$ and $v_2$ of $e_1$ and $e_2$ are in $U$ and are $U$-connected with each other (take the $\mathcal{T}$-path between  $v_1$ and $v_2$). Thus, we are in case $(\beta)$. Conversely, if we are in case $(\beta)$ during the algorithm, the current vertex $v$ is visited for the first time (or one of the vertices $v_1,~v_2$ would already be in $U$). Moreover, by the preceding point, one of the edges ($e_1$ or $e_2$) incident to $v$ is not in $\mathcal{T}$ and joins $v$ to one of its descendants. Hence, the current vertex $v$ is in $V_\beta$. \vspace{.2cm}\\
\ite During the DFS algorithm we have to make $n$ binary choices that will affect the outcome of the algorithm (case $(\beta)$). The other choices (case $(\gamma)$) do not affect the outcome of the algorithm. Therefore, there are $2^n$ possible outcomes.
\findem

\section{\textbf{Applications, extensions and open problems}} \label{section:perspectives}
\subsection{Random generation of triangulations} ~\\
The random generation of excursions of length $3n$ (with uniform distribution) reduces to the random generation of 1-dimensional walks of length $3n$ with steps +2, -1 starting and ending at 0 and remaining non-negative. The random generation of these walks is known to be feasible in linear time. (One just needs to generate a word of length $3n+1$ containing $n$ letters $c$ and $2n+1$ letters $\alpha$ and to apply the cycle lemma.) Given an excursion $w$, the construction of the 2-near-cubic marked-depth-map $\Phi(w)$ can be performed in linear time. Therefore, we have a linear time algorithm for the random generation (with uniform distribution) of bridgeless 2-near-cubic marked-depth-maps.  For any bridgeless 2-near-cubic map there are $2^n$ depth trees and $(n+1)$ possible marking. Therefore, if we drop the marking and the depth tree at the end of the process, we obtain a uniform distribution on bridgeless 2-near-cubic maps. This allows us to generate uniformly bridgeless cubic maps or, dually, loopless triangulations, in linear time.\\

\subsection{Kreweras walks ending at $(i,0)$ and $(i+2)$-near-cubic maps} ~\\
The Kreweras walks ending at  $(i,0)$ are the words $w$ on the alphabet $\{a,b,c\}$ with $|w|_a+i=|w|_b=|w|_c $ such that any suffix $w'$ of $w$ satisfies $|w'|_a +i \geq |w'|_c$ and $|w'|_b \geq |w'|_c$. There is a  very nice formula \cite{Kreweras:walks} giving the number of Kreweras walks of size $n$ (length $3n+2i$) ending at $(i,0)$:
\begin{eqnarray}\label{eq:kni}
k_{n,i}&=& \frac{4^n(2i+1)}{(n+i+1)(2n+2i+1)}{2i \choose i}{3n+2i \choose n}.
\end{eqnarray}
There is also a similar formula \cite{Mullin:triangulation-nonsep} for non-separable $(i+2)$-near-cubic maps of size $n$ ($3n+2i+1$ edges):
\begin{eqnarray}\label{eq:tni}
c_{n,i}&=& \frac{2^n(2i+1)}{(n+i+1)(2n+2i+1)}{2i \choose i}{3n+2i \choose n}.
\end{eqnarray}

In this subsection, we show that the bijection $\Phi$ 
(Definition \ref{def:Phi-prime}) can be extended to Kreweras walks ending at  $(i,0)$. This gives a bijective correspondence explaining why $k_{n,i}=2^nc_{n,i}$. \\

\begin{figure}[h!]
\begin{center}
\input{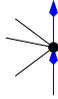}
\caption{The tree-growing map $M_\bullet^i$ when $i=3$.}\label{fig:i-seed}
\end{center}
\end{figure}

Consider the tree-growing map $M_\bullet^i$ reduced to a vertex, a root, a head and $i$ left legs (Figure \ref{fig:i-seed}). We define the image of a  Kreweras walk $w$ ending at  $(i,0)$  as the map obtained by closing $\varphi_w(M_\bullet^i)$. We get the following extension of Theorem~\ref{thm:bijection}.

\begin{thm}\label{thm:bijection-again}
The mapping $\Phi$ is a bijection between Kreweras walks of size $n$ (length $3n+2i$) ending at $(i,0)$  and  non-separable $(i+2)$-near-cubic maps of size $n$  ($3n+2i+1$ edges) marked on the root-edge with a depth tree that contains the edge following the root in counterclockwise order around the root-vertex.
\end{thm}

By Theorem \ref{thm:counting-depth-tree}, there are $2^n$ such depth trees. Consequently, we obtain the following corollary:

\begin{cor}
The number $k_{n,i}$ of Kreweras walks of size $n$ ending at $(i,0)$ and the number $c_{n,i}$ of  non-separable $(i+2)$-near-cubic maps of size~$n$ are related by the equation $k_{n,i}=2^n c_{n,i}$.
\end{cor}

One can define the counterpart of excursions for Kreweras walks ending at $(i,0)$. These are the walks obtained when one chooses an external edge in a non-separable $(i+2)$-near-cubic depth-map such that the edge following the root is in the tree and applies the mapping $\Psi=\Phi^{-1}$. Alas, we have found no simple characterization of this set of walks nor any bijective proof explaining why this set has cardinality $\displaystyle \frac{4^n(2i+1)}{(2n+2i+1)}{2i \choose i}{3n+2i \choose n}$.\vspace{.5cm}\\

\noindent \textbf{Acknowledgments:} I thank Mireille Bousquet-Mélou for steady encouragements and support.

\bibliography{../../../biblio/allref}
\bibliographystyle{plain}













\end{document}